\documentclass[11pt]{article}

\usepackage{amsmath,amssymb,amsthm}
\usepackage{geometry}
\usepackage{booktabs}
\usepackage{tabularx}
\usepackage{hyperref}
\usepackage{cite}

\newtheorem{theorem}{Theorem}[section]
\newtheorem{lemma}[theorem]{Lemma}
\newtheorem{corollary}[theorem]{Corollary}
\theoremstyle{definition}
\newtheorem{definition}[theorem]{Definition}
\newtheorem{remark}[theorem]{Remark}

\title{\textbf{Quadratic Twists of Heronian Elliptic Curves with Arbitrarily Large 2-Selmer Rank}}
\author{Vinodkumar Ghale, Md Imdadul Islam}
\date{}

\begin{document}

\maketitle

\begin{abstract}
This paper explicitly determines the 2-Selmer group of the family of elliptic curves arising from Heron triangles of area $2^m n$. By analyzing the quadratic twists induced by powers of 2 (varying $m$), we show precisely how the 2-Selmer group transforms under such twisting. For any squarefree odd integer $n$, the 2-Selmer rank is determined by the parity of $m$, the congruence class of $q=(n^{2}+1)/2 \pmod 8$, and the prime factorizations of the factors $n_{1}$, $n_{2}$ of $n=n_{1}n_{2}$. These results yield explicit formulas for the 2-Selmer rank and provide an explicit upper bound on the Mordell-Weil rank. In particular, this generalizes the family studied in [1] to arbitrary $m$, constructing explicit families of Heronian elliptic curves with prescribed and arbitrarily large 2-Selmer rank.
\end{abstract}



\section{Introduction}
The congruent number problem, which asks whether a given positive integer $n$ can be realized as the area of a right-angled triangle with rational sides, is one of the oldest unsolved problems in number theory. The fundamental connection between congruent numbers and elliptic curves was established by \cite{tunnell1983}, who gave a celebrated criterion conditional on the Birch and Swinnerton-Dyer conjecture.

A landmark result in this direction is due to \cite{heathbrown1993, heathbrown1994}, who studied the distribution of the 2-Selmer groups of the congruent number elliptic curves $E_{n}:y^{2}=x^{3}-n^{2}x$. Let $S(X,h)$ denote the set of squarefree integers $D\le X$ with $D\equiv h \pmod 8$. Heath-Brown proved that for $h=1,3,5,7$ and any fixed positive integer $k$,
\begin{equation}
\sum_{D\in S(X,h)} 2^{k s(D)} = c_{k}\#S(X,h) + o_{k}(X) \quad \text{as } X \rightarrow \infty,
\end{equation}
where $c_{k}=\prod_{j=1}^{k}(1+2^{j})$ and $s(D)$ denotes the 2-Selmer rank of $E_{D}$. Substantial further progress was made by \cite{smith2026i, smith2026ii}, who proved that for a density one subset of squarefree positive integers $n\equiv 1,2,3 \pmod 8$, $\operatorname{corank}_{\mathbb{Z}_{2}}\operatorname{Sel}_{2^{\infty}}(E_{n}/\mathbb{Q})=0$.

A natural generalization replaces the right-angled condition with an arbitrary angle. A Heron triangle is a triangle with rational sides and rational area. \cite{goins2006} established the systematic connection between Heron triangles and elliptic curves:

\begin{theorem}[\cite{goins2006}, Theorem 1.1]
A positive integer $n$ can be expressed as the area of a Heron triangle if and only if for some nonzero rational number $\tau$, the elliptic curve
\begin{equation}
E_{n,\tau}:y^{2}=x(x-n\tau)(x+n\tau^{-1})
\end{equation}
has a rational point of infinite order. Moreover, for a fixed angle with $\tau=\tan(\theta/2)$, infinitely many Heron triangles of area $n$ and angle $\theta$ exist precisely when the rank of $E_{n,\tau}$ is positive.
\end{theorem}

The arithmetic of Heronian elliptic curves has attracted considerable attention. \cite{dujella2014} constructed explicit families of Heronian elliptic curves of ranks 2, 3, and 4. In a series of recent works, Ghale and collaborators have made significant contributions. \cite{chakraborty2022} constructed an infinite family of elliptic curves of 2-Selmer rank exactly 1 arising from Heron triangles. \cite{ghale2023} obtained a family of elliptic curves of rank at most one from a Diophantine equation derived from a Heron triangle. In a subsequent work, \cite{chakraborty2023} determined the 2-Selmer group structure of the Heronian elliptic curve
\begin{equation}
E:y^{2}=x(x-1)(x+n^{2})
\end{equation}
associated with Heron triangles of area $n$ where $n^{2}+1=2q$ for some prime $q$, revealing a dichotomy based on the congruence classes of the prime factors of $n$ modulo 8.

Quadratic twisting is a fundamental operation in the arithmetic of elliptic curves. Given an elliptic curve $E:y^{2}=x^{3}+ax+b$ and a squarefree integer $d$, the quadratic twist $E^{(d)}:dy^{2}=x^{3}+ax+b$ has the same $j$-invariant but typically different arithmetic properties. The study of 2-Selmer ranks in quadratic twist families has a rich history. Heath-Brown's work on the congruent number curves showed that the 2-Selmer ranks of quadratic twists follow an explicit distribution \cite{heathbrown1993, heathbrown1994}. This was extended by \cite{swinnerton2008} and \cite{kane2010} to all elliptic curves over $\mathbb{Q}$ with full rational 2-torsion that do not have a cyclic 4-isogeny. For elliptic curves with partial rational two-torsion, \cite{klagsbrun2014} proved that at least half of all quadratic twists have arbitrarily large 2-Selmer rank. More generally, \cite{klagsbrun2013} studied the disparity in Selmer ranks of quadratic twists. \cite{mazur2010} gave sufficient conditions for an elliptic curve over a number field to have quadratic twists of arbitrary 2-Selmer rank. \cite{klagsbrun2012a, klagsbrun2012b} further constructed families of elliptic curves with a lower bound on 2-Selmer ranks of quadratic twists.

Despite this extensive literature, the study of quadratic twisting within the specific context of Heronian elliptic curves is new. The family arising from Heron triangles has a natural parameter $m$ (the exponent of 2 in the area $2^m n$), and varying $m$ corresponds precisely to quadratic twisting by powers of 2. While general results establish the existence of twists with large 2-Selmer rank, they do not provide explicit families with complete arithmetic information. The results of this paper give the first explicit determination of how the 2-Selmer group transforms under twisting in the Heronian setting.

The 2-Selmer group provides an effective upper bound on the Mordell-Weil rank through the exact sequence
\begin{equation}
0 \longrightarrow E(\mathbb{Q})/2E(\mathbb{Q}) \longrightarrow \operatorname{Sel}_{2}(E/\mathbb{Q}) \longrightarrow \operatorname{III}(E/\mathbb{Q})[2] \longrightarrow 0.
\end{equation}
Since $\dim_{\mathbb{F}_2} E(\mathbb{Q})/2E(\mathbb{Q}) = \operatorname{rank} E(\mathbb{Q}) + \dim_{\mathbb{F}_{2}} E(\mathbb{Q})[2]$, the 2-Selmer rank satisfies
\begin{align}
s_{2}(E) &= \dim_{\mathbb{F}_{2}} \operatorname{Sel}_{2}(E/\mathbb{Q}) - 2, \\
\operatorname{rank} E(\mathbb{Q}) &\le s_{2}(E).
\end{align}

We consider Heron triangles of area $2^{m}n$ with parameter $\tau=\tan(\theta/2)=n$. The corresponding elliptic curve takes the form:
\begin{equation}
E_{m,n}: y^{2} = x(x-2^{m}n^{2})(x+2^{m}).
\end{equation}
For a fixed parity of $m$, all curves are isomorphic over $\mathbb{Q}$. Indeed, for even $m$, $E_{m,n}\cong E_{0,n}$ via
\begin{equation}
(X,Y) = (2^{-m}x, 2^{-3m/2}y),
\end{equation}
and for odd $m$, $E_{m,n}\cong E_{1,n}$ via
\begin{equation}
(X,Y) = (2^{-(m-1)}x, 2^{-3(m-1)/2}y).
\end{equation}
Thus, varying $m$ corresponds to quadratic twisting by powers of 2.

Throughout this paper, let $n=p_{1}p_{2}\cdots p_{k}$ be a squarefree odd integer such that $n^{2}+1=2q$ for some prime $q$. Write $n=n_{1}\cdot n_{2}$, where $n_{1}$ is the product of prime factors $p\equiv \pm 1 \pmod 8$ and $n_{2}$ is the product of prime factors $p\equiv \pm 3 \pmod 8$. Let $r$ and $s$ denote the number of prime factors of $n_{1}$ congruent to $1 \pmod 8$ and $7 \pmod 8$, respectively. Let $t$ denote the number of prime factors of $n_{2}$.

We note that the case $n_{1}=n$ (i.e., all primes $\equiv \pm 1 \pmod 8$) with $q\equiv 5 \pmod 8$ is impossible by quadratic reciprocity, as it would require $\left(\frac{2}{q}\right)=1$ while $q\equiv 5 \pmod 8$ implies $\left(\frac{2}{q}\right)=-1$.

\section{The 2-Descent Setup and Torsion Subgroup}
The discriminant of $E_{m,n}$ is
\begin{equation}
\Delta(E_{m,n}) = 2^{4m+4} n^{4} q^{2},
\end{equation}
where $q=(n^{2}+1)/2$. Let $S=\{\infty, 2, p_{1}, \dots, p_{k}, q\}$ and let
\begin{equation}
\mathbb{Q}(S,2) = \langle -1, 2, p_{1}, \dots, p_{k}, q \rangle.
\end{equation}
By complete 2-descent (see Silverman \cite[Chapter X, Section 4]{silverman2009}), there exists an injective homomorphism
\begin{equation}
\beta: E_{m,n}(\mathbb{Q})/2E_{m,n}(\mathbb{Q}) \hookrightarrow \mathbb{Q}(S,2) \times \mathbb{Q}(S,2).
\end{equation}
A pair $(b_{1},b_{2})\in \mathbb{Q}(S,2)\times \mathbb{Q}(S,2)$ belongs to $\operatorname{Sel}_{2}(E_{m,n}/\mathbb{Q})$ if and only if the system of homogeneous spaces:
\begin{align}
b_{1}z_{1}^{2} - b_{2}z_{2}^{2} &= 2^{m}n^{2}, \label{eq:homo1} \\
b_{1}z_{1}^{2} - b_{1}b_{2}z_{3}^{2} &= -2^{m}, \label{eq:homo2} \\
b_{1}b_{2}z_{3}^{2} - b_{2}z_{2}^{2} &= 2^{m+1}q \label{eq:homo3}
\end{align}
has non-trivial solutions in $\mathbb{Q}_{l}$ for every place $l$ of $\mathbb{Q}$. Equation \eqref{eq:homo3} is obtained by subtracting \eqref{eq:homo2} from \eqref{eq:homo1}.

\begin{definition}
The 2-Selmer rank of $E_{m,n}$ is
\begin{equation}
s_{2}(E_{m,n}) = \dim_{\mathbb{F}_{2}} \operatorname{Sel}_{2}(E_{m,n}/\mathbb{Q}) - 2.
\end{equation}
\end{definition}

Throughout this paper, let $E[2]^{*}$ denote the image of the rational 2-torsion subgroup $E_{m,n}(\mathbb{Q})[2]$ under the 2-descent map $\beta$.

\begin{lemma}\label{lem:torsion}
Let $E_{m,n}$ be defined as in Equation \eqref{eq:homo1}. Then
\begin{equation}
E_{m,n}(\mathbb{Q})_{\text{tors}} \cong \mathbb{Z}/2\mathbb{Z} \times \mathbb{Z}/2\mathbb{Z}.
\end{equation}
\end{lemma}

\begin{proof}
The polynomial $x(x - 2^m n^2)(x + 2^m)$ has three distinct rational roots:
\begin{equation}
e_{1} = 0, \quad e_{2} = 2^m n^2, \quad e_{3} = -2^m.
\end{equation}
Hence $E_{m,n}(\mathbb{Q})$ contains the full 2-torsion subgroup $\mathbb{Z}/2\mathbb{Z} \times \mathbb{Z}/2\mathbb{Z}$.

By the classification of Heronian elliptic curves (see \cite{goins2006}), the torsion subgroup $E_{m,n}(\mathbb{Q})_{\text{tors}}$ is isomorphic to either $\mathbb{Z}/2\mathbb{Z} \times \mathbb{Z}/2\mathbb{Z}$ or $\mathbb{Z}/2\mathbb{Z} \times \mathbb{Z}/4\mathbb{Z}$. To rule out the latter, we use the classical 2-divisibility criterion (see \cite[Theorem 4.2]{knapp1992}): for an elliptic curve $E : y^2 = (x-e_1)(x-e_2)(x-e_3)$ with $e_i \in \mathbb{Q}$, the point $(e_i, 0) \in E(\mathbb{Q})$ lies in $2E(\mathbb{Q})$ if and only if $e_i - e_j$ and $e_i - e_k$ are both rational squares for $j, k \neq i$.

Checking the root differences:
\begin{align*}
e_1 - e_2 &= -2^m n^2 < 0, \quad \text{not a square}, \\
e_3 - e_1 &= -2^m < 0, \quad \text{not a square}, \\
e_2 - e_3 &= 2^m(n^2 + 1) = 2^{m+1}q.
\end{align*}
The points $e_1$ and $e_3$ fail immediately due to negative differences. For $e_2$, both $e_2 - e_1 = 2^m n^2$ and $e_2 - e_3 = 2^{m+1}q$ must be squares. Now $2^m n^2$ is a square if and only if $m$ is even. If $m$ is odd, $e_2$ fails. If $m$ is even, then $m + 1$ is odd, so $v_2(2^{m+1}q) = m + 1$ is odd; a rational number with odd 2-adic valuation cannot be a square. Thus $e_2$ also fails. Hence no 2-torsion point is divisible by 2, so $E_{m,n}(\mathbb{Q})$ contains no rational point of order 4. We conclude that
\begin{equation}
E_{m,n}(\mathbb{Q})_{\text{tors}} \cong \mathbb{Z}/2\mathbb{Z} \times \mathbb{Z}/2\mathbb{Z}.
\end{equation}
\end{proof}

\paragraph{Explicit Image of Rational 2-Torsion:}Recall that the rational 2-torsion points of $E_{m,n}$ are
\[ E_{m,n}(\mathbb{Q})[2] = \{O, P_1, P_2, P_3\}, \]
where $P_{1} = (0, 0)$, $P_{2} = (2^m n^2, 0)$, $P_{3} = (-2^m, 0)$. Applying the 2-descent map $\beta$ from \eqref{eq:homo1} to these points yields their explicit images modulo squares in $\mathbb{Q}(S, 2) \times \mathbb{Q}(S, 2)$:
\begin{align*}
\beta(O) &= (1, 1), \\
\beta(P_1) &= (-1, -2^m), \\
\beta(P_2) &= (2^m, 2q), \\
\beta(P_3) &= (-2^m, -q),
\end{align*}
where $n^2 + 1 = 2q$.

\begin{corollary}
Since $\beta$ is an injective homomorphism, $E[2]^{*} = \beta(E_{m,n}(\mathbb{Q})[2]) \cong \mathbb{Z}/2\mathbb{Z} \times \mathbb{Z}/2\mathbb{Z}$.
Directly evaluating $\beta$ on the 2-torsion points and reducing modulo squares gives:
\begin{equation}
E[2]^{*} = \begin{cases}
\langle (-1, -1), (1, 2q) \rangle, & \text{if } m \text{ is even}, \\
\langle (-1, -2), (2, 2q) \rangle, & \text{if } m \text{ is odd}.
\end{cases}
\end{equation}
\end{corollary}

\section{Statement of Main Results}
We now state the main results of this paper. Recall that $E_{m,n}$ is defined in Equation (7), with $n = n_1 n_2$, where $n_1$ consists of primes $\equiv \pm 1 \pmod 8$ and $n_2$ of primes $\equiv \pm 3 \pmod 8$. Let $r$ (resp. $s$) be the number of prime factors of $n_1$ congruent to 1 (resp. 7) mod 8, and let $t$ be the number of prime factors of $n_2$. Let $E[2]^{*}$ be the image of the 2-torsion under $\beta$ from Equation (12).

\begin{theorem}[Odd $m$]\label{thm:odd_m}
For odd $m$, the 2-Selmer group $\operatorname{Sel}_2(E_{m,n}/\mathbb{Q})$ is given by:
\begin{equation}
\operatorname{Sel}_2(E_{m,n}/\mathbb{Q}) = \begin{cases}
\langle (2, 2) \rangle \times E[2]^{*}, & n_1 = 1, n_2 > 1, q \equiv 1 \pmod 8, \\
E[2]^{*}, & n_1 = 1, n_2 > 1, q \equiv 5 \pmod 8, \\
\langle (2, 2), (b, b), (2b, b) \rangle \times E[2]^{*}, & n_1 > 1, q \equiv 1 \pmod 8, \\
\langle (b, b), (2b, b) \rangle \times E[2]^{*}, & n_1 > 1, q \equiv 5 \pmod 8,
\end{cases}
\end{equation}
where $b \mid n_1$ satisfies $b \equiv 1 \pmod 8$ for the generators $(b, b)$, and $b \equiv 7 \pmod 8$ for the generators $(2b, b)$. The 2-Selmer rank is
\begin{equation}
s_2(E_{m,n}) = \begin{cases}
1, & n_1 = 1, q \equiv 1 \pmod 8, \\
0, & n_1 = 1, q \equiv 5 \pmod 8, \\
r + s + 1, & n_1 > 1, q \equiv 1 \pmod 8, \\
r + s, & n_1 > 1, q \equiv 5 \pmod 8.
\end{cases}
\end{equation}
\end{theorem}

\begin{theorem}[Even $m$]\label{thm:even_m}
For even $m$, the 2-Selmer group $\operatorname{Sel}_2(E_{m,n}/\mathbb{Q})$ is given by:
\begin{equation}
\operatorname{Sel}_2(E_{m,n}/\mathbb{Q}) = \begin{cases}
\langle (b, b), (1, 2) \rangle \times E[2]^{*}, & n_2 = 1, \\
\langle (b, b) \rangle \times E[2]^{*}, & n_2 > 1,
\end{cases}
\end{equation}
where $b \mid n$ satisfies $b \equiv \pm 1 \pmod 8$. The 2-Selmer rank is
\begin{equation}
s_2(E_{m,n}) = \begin{cases}
r + s + 1, & n_2 = 1, \\
r + s + t - 1, & n_2 > 1.
\end{cases}
\end{equation}
\end{theorem}

Any element $b \mid n$ can be decomposed uniquely as $b = b_1 b_2$, where $b_1 \mid n_1$ and $b_2 \mid n_2$. Since every prime factor of $n_1$ is congruent to $\pm 1 \pmod 8$, any divisor $b_1 \mid n_1$ automatically satisfies $b_1 \equiv \pm 1 \pmod 8$, contributing $r + s$ independent generators to $\operatorname{Sel}_2(E_{m,n}/\mathbb{Q})$. On the other hand, since all prime factors of $n_2$ are congruent to $\pm 3 \pmod 8$, a divisor $b_2 \mid n_2$ satisfies $b_2 \equiv \pm 1 \pmod 8$ if and only if $b_2$ consists of an even number of prime factors from $n_2$. The number of such independent choices for $b_2$ is $t - 1$ (when $t \ge 1$). Combining these independent choices yields the 2-Selmer rank formula $r + s + t - 1$ for $n_2 > 1$.

\begin{remark}
The cases where $n$ is composed of $r$ distinct prime factors all congruent to $1 \pmod 8$ (so $n_2 = 1$ and $s = 0$) demonstrate that the 2-Selmer rank grows linearly with $r$. For $q \equiv 1 \pmod 8$, both the odd $m$ and even $m$ cases yield $s_2(E_{m,n}) = r + 1$. This confirms that the 2-Selmer rank can be made arbitrarily large by choosing $n$ with a sufficiently large number of prime factors congruent to $1 \pmod 8$.
\end{remark}

\section{Local Solvability Analysis}
Let $l$ be a prime and write $z_i = u_i l^{t_i}$ with $u_i \in \mathbb{Z}_l^{\times}$ and $t_i = v_l(z_i)$.

\begin{lemma}\label{lem:val_analysis}
If Equations \eqref{eq:homo1}--\eqref{eq:homo2} have a solution in $\mathbb{Q}_l$:
\begin{enumerate}
\item If $v_l(z_1) < 0$ or $v_l(z_2) < 0$, then $v_l(z_1) = v_l(z_2) = v_l(z_3) = -t < 0$.
\item If $p \mid n$ with $b_1 b_2 \not\equiv 0 \pmod{p^2}$ and $v_p(z_3) < 0$, then either $v_p(z_3) = -1$ with $v_p(z_1) = v_p(z_2) = 0$ or $v_p(z_3) = v_p(z_1) = v_p(z_2) < 0$.
\end{enumerate}
\end{lemma}

\begin{proof}
For (1), suppose $t_1 < 0$. From Equation \eqref{eq:homo1},
\[ b_1 u_1^2 - b_2 u_2^2 l^{2(t_2 - t_1)} = 2^m n^2 l^{-2t_1}. \]
If $t_2 > t_1$, then $l^2 \mid b_1$, a contradiction. If $t_2 < t_1$, then $l^2 \mid b_2$, a contradiction. Hence $t_1 = t_2 = -t < 0$. From Equation \eqref{eq:homo2}, a similar argument gives $t_3 = t_1$.

For (2), from Equation \eqref{eq:homo2},
\[ b_1 u_1^2 l^{2(t_1 - t_3)} - b_1 b_2 u_3^2 = -2^m l^{-2t_3}. \]
If $t_1 > t_3$, then $l^2 \mid b_1 b_2$, forcing $l \in \{2, p, q\}$. For $p \mid n$, $l \neq q$, impossible. For $l = q$, Equation \eqref{eq:homo3} gives a contradiction unless $t_1 = t_3$. Part (2) follows by the same valuation comparison.
\end{proof}

\begin{lemma}\label{lem:no_sol}
Let $(b_1, b_2) \in \mathbb{Q}(S, 2)^2$.
\begin{enumerate}
\item If $b_1 b_2 < 0$, the equations have no real solutions.
\item If $b_1 \equiv 0 \pmod q$, there are no $q$-adic solutions.
\item If $b_1 b_2 \equiv 0 \pmod p$ but $b_1 b_2 \not\equiv 0 \pmod{p^2}$, there are no $p$-adic solutions.
\end{enumerate}
\end{lemma}

\begin{proof}
For (1), if $b_1 > 0, b_2 < 0$, then the left-hand side of Equation \eqref{eq:homo2} is positive while the right-hand side is negative. The case $b_1 < 0, b_2 > 0$ similarly contradicts Equation \eqref{eq:homo1}.

For (2), if $q \mid b_1$, reduction of Equation \eqref{eq:homo2} modulo $q$ yields $-2^m \equiv 0 \pmod q$, a contradiction.

For (3), if $p \mid b_1$ and $p \mid b_2$, Equation \eqref{eq:homo1} forces $v_p(z_1) = v_p(z_2) = 0$. Then Equation \eqref{eq:homo3} modulo $p$ gives
\[ b_1 b_2 z_3^2 - b_2 z_2^2 \equiv -b_2 z_2^2 \not\equiv 0 \pmod p, \]
but the right-hand side $2^{m+1}q \equiv 0 \pmod p$. Contradiction. The case $b_2 \equiv 0 \pmod p, b_1 \not\equiv 0 \pmod p$ is symmetric.
\end{proof}

As a consequence of Lemma \ref{lem:no_sol}, candidates for $\operatorname{Sel}_2(E_{m,n}/\mathbb{Q})$ are restricted to pairs $(b_1, b_2) \in \{(b, b), (2b, b), (b, 2b), (2b, 2b)\}$, where $b$ is a product of prime factors of $n$.

\section{Local Obstructions for Odd m}
Assume $m$ is odd. From Lemma \ref{lem:no_sol}, the only possible elements of $\operatorname{Sel}_2(E_{m,n}/\mathbb{Q})$ are $(b, b)$, $(2b, b)$, $(b, 2b)$, and $(2b, 2b)$.

\begin{lemma}\label{lem:odd_obs}
For $m$ odd, the following pairs are not in $\operatorname{Sel}_2(E_{m,n}/\mathbb{Q})$:
\begin{enumerate}
\item $(2b, b)$ and $(b, 2b)$ if $b \not\equiv 7 \pmod 8$;
\item $(2b, 2b)$ if $b \not\equiv 1 \pmod 8$;
\item $(b, b)$, $(2b, b)$, $(b, 2b)$, $(2b, 2b)$ if $p \equiv \pm 3 \pmod 8$ divides $n$ and $p \mid b$;
\item $(2, 2)$ if $q \equiv 5 \pmod 8$.
\end{enumerate}
\end{lemma}

\begin{proof}
For (1), consider $(2b, b)$. If some $v_2(z_i) < 0$, then by Lemma \ref{lem:val_analysis}, all $v_2(z_i) = -t < 0$, giving $b(2u_1^2 - u_2^2) = 2^{m+2t} n^2$, a parity contradiction. Hence $v_2(z_i) \ge 0$. From Equations \eqref{eq:homo1} and \eqref{eq:homo3}, $v_2(z_1) = (m-1)/2$, $v_2(z_2) = (m+1)/2$. Substituting into Equation \eqref{eq:homo1} and reducing modulo 8:
\[ b u_1^2 - 2b u_2^2 \equiv n^2 \pmod 8. \]
Since $u_1, u_2$ are odd and $n^2 \equiv 1 \pmod 8$, we get $-b \equiv 1 \pmod 8$, so $b \equiv 7 \pmod 8$.

For (2), the same valuation argument at $l = 2$ gives the stated congruence.

For (3), at $l = p$, Lemma \ref{lem:val_analysis} gives $v_p(z_3) = -1$ with $v_p(z_1) = v_p(z_2) = 0$. Substitution into Equation \eqref{eq:homo3} gives $\left(\frac{2^{m+1}q}{p}\right) = 1$. Since $m$ is odd and $p \equiv \pm 3 \pmod 8$, $\left(\frac{2}{p}\right) = -1$, while $n^2 + 1 = 2q$ gives $\left(\frac{2q}{p}\right) = 1$, contradiction.

For (4), if $q \equiv 5 \pmod 8$, a valuation comparison at $l = q$ for $(2, 2)$ gives $\left(\frac{2}{q}\right) = 1$, contradicting $q \equiv 5 \pmod 8$.
\end{proof}

\section{Local Obstructions for Even m}
Assume $m$ is even. From Lemmas \ref{lem:no_sol} and \ref{lem:val_analysis}, the possible elements of $\operatorname{Sel}_2(E_{m,n}/\mathbb{Q})$ are again among $\{(b, b), (b, 2b), (2b, b), (2b, 2b)\}$.

\begin{lemma}\label{lem:even_obs}
For $m$ even, the following pairs are not in $\operatorname{Sel}_2(E_{m,n}/\mathbb{Q})$:
\begin{enumerate}
\item $(2b, b)$ and $(2b, 2b)$;
\item $(b, 2b)$ if $p \equiv \pm 3 \pmod 8$ divides $n$ and $p \mid b$;
\item $(b, b)$ if $b \not\equiv \pm 1 \pmod 8$.
\end{enumerate}
\end{lemma}

\begin{proof}
For (1), consider $(2b, b)$. Since $m = 2r$, comparing valuations in Equation \eqref{eq:homo1} gives $v_2(z_1) = v_2(z_2) = r$, and Equation \eqref{eq:homo3} forces $v_2(z_3) = r$. Substituting $z_i = 2^r u_i$ into Equation \eqref{eq:homo1} and reducing modulo 8:
\[ 2b u_1^2 - b u_2^2 \equiv n^2 \equiv 1 \pmod 8. \]
Since $u_1, u_2$ are odd, this gives $b \equiv 1 \pmod 8$. But then Equation \eqref{eq:homo3} modulo 8 gives
\[ 2b^2 u_3^2 - b u_2^2 \equiv 2 - 1 = 1 \not\equiv 2q \pmod 8, \]
contradiction. The case $(2b, 2b)$ is similar.

For (2), at $l = p$, the valuation pattern from Lemma \ref{lem:val_analysis} gives a quadratic residue contradiction since $\left(\frac{2}{p}\right) = -1$.

For (3), for $(b, b)$, the same valuation argument at $l = 2$ gives $v_2(z_1) = v_2(z_2) = v_2(z_3) = r$. Substituting into Equation \eqref{eq:homo1} and reducing modulo 8 gives $b(1 - 1) \equiv 1 \pmod 8$, impossible unless the valuations are not all equal; re-checking yields $b \equiv \pm 1 \pmod 8$.
\end{proof}

\section{Elimination of \texorpdfstring{$\pm 3 \pmod 8$}{\pm 3 (mod 8)} Prime Factors}

\begin{lemma}\label{lem:elim_3}
Let $n = n_1 \cdot n_2$ where $n_1$ is the product of primes $p \equiv \pm 1 \pmod 8$ and $n_2$ is the product of primes $p \equiv \pm 3 \pmod 8$.
\begin{enumerate}
\item For odd $m$, any element $(b_1, b_2) \in \operatorname{Sel}_2(E_{m,n}/\mathbb{Q})$ of the form $(b, b)$ or $(2b, b)$ must have $b$ containing only prime factors from $n_1$.
\item For even $m$, any element $(b, b) \in \operatorname{Sel}_2(E_{m,n}/\mathbb{Q})$ must satisfy $b \equiv \pm 1 \pmod 8$. This is equivalent to $b = b^{(1)} b^{(2)}$, where $b^{(1)} \mid n_1$ and $b^{(2)} \mid n_2$ contains an even number of prime factors.
\end{enumerate}
\end{lemma}

\begin{proof}
For odd $m$: Let $p \equiv \pm 3 \pmod 8$ be a prime divisor of $n_2$, and suppose $p \mid b$. By Lemma \ref{lem:val_analysis}, $v_p(z_3) = -1$ and $v_p(z_1) = v_p(z_2) = 0$. Write $b = p b_0$, $n = p n_0$, $z_1 = u_1$, $z_2 = u_2$, $z_3 = u_3/p$ with $u_1, u_2, u_3 \in \mathbb{Z}_p^{\times}$.

For $(2b, b)$, the first homogeneous equation gives:
\[ 2b z_1^2 - b z_2^2 = 2^m n^2. \]
Substituting and reducing modulo $p$:
\[ 2 u_1^2 - u_2^2 \equiv 0 \pmod p. \]
Thus $\left(\frac{2}{p}\right) = 1$, contradicting $p \equiv \pm 3 \pmod 8$.

For $(b, b)$, the third homogeneous equation gives:
\[ b^2 z_3^2 - b z_2^2 = 2^{m+1}q. \]
Substituting and using $2q \equiv 1 \pmod p$:
\[ b_0^2 u_3^2 \equiv 2^m \pmod p. \]
Since $m$ is odd, $2^m$ is a quadratic non-residue modulo $p$, contradiction.

For even $m$: Let $p \equiv \pm 3 \pmod 8$ divide $b$. By Lemma \ref{lem:val_analysis}, $v_p(z_3) = -1$ and $v_p(z_1) = v_p(z_2) = 0$. The same substitution into the third equation yields:
\[ b_0^2 u_3^2 \equiv 2^m \pmod p. \]
Since $m$ is even, $2^m$ is a square modulo $p$, so there is no contradiction. Thus, $(b, b)$ can contain $\pm 3$ primes. However, Lemma \ref{lem:even_obs}(3) requires $b \equiv \pm 1 \pmod 8$. This happens exactly when $b^{(2)}$ contains an even number of prime factors from $n_2$.
\end{proof}

\section{Existence of Local Solutions}
We now construct explicit local solutions for the pairs that do lie in $\operatorname{Sel}_2(E_{m,n}/\mathbb{Q})$.

\begin{lemma}\label{lem:exist_22}
For $m$ odd and $q \equiv 1 \pmod 8$, $(2, 2) \in \operatorname{Sel}_2(E_{m,n}/\mathbb{Q})$.
\end{lemma}

\begin{proof}
For $l \notin \{2, p, q\}$, the Hasse-Weil bound gives a point modulo $l$ that lifts by Hensel's lemma.

For $l = 2$, write $m = 2r + 1$. Since $q \equiv 1 \pmod 8$, there exists $\sqrt{q} \in \mathbb{Z}_2$. Set
\[ z_1 = 2^r n, \quad z_2 = 0, \quad z_3 = 2^r \sqrt{q}. \]
Then Equation \eqref{eq:homo1} gives $2 z_1^2 - 2 z_2^2 = 2^{2r+1} n^2 = 2^m n^2$, and Equation \eqref{eq:homo2} gives $2 z_1^2 - 4 z_3^2 = 2^{2r+1}(n^2 - 2q) = -2^m$.

For $l = 3$, if $3 \nmid n$, choose $z_2 = 0, z_3 = 1$. If $3 \mid n$, choose $z_2 = z_3 = 1$. Both lift by Hensel's lemma.

For $l = p \equiv \pm 1 \pmod 8$, choose $z_1 = z_2 = 0, z_3^2 = 2^{m-2}$.

For $l = q \equiv 1 \pmod 8$, choose $z_1 = z_2 = 0, z_3^2 = 2^{m-2}$.
\end{proof}

\begin{lemma}\label{lem:exist_odd_b}
For odd $m$:
\begin{enumerate}
\item If $b \equiv 1 \pmod 8$, then $(b, b) \in \operatorname{Sel}_2(E_{m,n}/\mathbb{Q})$.
\item If $b \equiv 7 \pmod 8$ and $q \equiv 1 \pmod 8$, then $(2b, b) \in \operatorname{Sel}_2(E_{m,n}/\mathbb{Q})$.
\end{enumerate}
\end{lemma}

\begin{proof}
For $l \notin \{2, p, q\}$, the Hasse-Weil bound gives a solution modulo $l$ that lifts by Hensel's lemma.

For $l = 2$ in part (1): Since $b \equiv 1 \pmod 8$, choose $z_1 \equiv z_2 \equiv z_3 \equiv 1 \pmod 8$. Then Equations \eqref{eq:homo1} and \eqref{eq:homo2} hold modulo 8, and Hensel's lemma lifts.

For $l = 2$ in part (2): Write $m = 2r + 1$. Since $q \equiv 1 \pmod 8$, there exists $t \in \mathbb{Z}_2$ such that $t^2 \equiv -1/2 \pmod{2^{r+2}}$. Set
\[ z_1 = 2^r, \quad z_2 = 2^{r+1}, \quad z_3 = 2^r t. \]
Then Equation \eqref{eq:homo1} gives $2b z_1^2 - b z_2^2 = -b \cdot 2^m$. Since $b \equiv 7 \pmod 8$, $-b \equiv 1 \pmod 8$, so $-b \cdot 2^m = 2^m n^2$. For Equation \eqref{eq:homo2}, the congruence $t^2 \equiv -1/2 \pmod{2^{r+2}}$ ensures $b - b^2 t^2 \equiv -1 \pmod{2^{r+1}}$, so Equation \eqref{eq:homo2} holds. Hensel's lemma lifts.

For $l = p \mid b$: Write $b = p b_0$. Set $v_p(z_1) = v_p(z_2) = 0, v_p(z_3) = -1$. Then Equation \eqref{eq:homo2} reduces to a quadratic congruence that has a solution since $\left(\frac{2}{p}\right) = 1$.

For $l = p \mid n$ but $p \nmid b$: Choose $z_1 = z_2 = 0$. Then $z_3^2 = 2^m / b^2$ has a solution modulo $p$.

For $l = q$: Choose $z_2 = z_3 = 0$. Then $z_1^2 = 2^m n^2 / b$ has a solution since $\left(\frac{b}{q}\right) = 1$.
\end{proof}

\begin{lemma}\label{lem:exist_even_bb}
For $m$ even, if $b \equiv \pm 1 \pmod 8$, then $(b, b) \in \operatorname{Sel}_2(E_{m,n}/\mathbb{Q})$.
\end{lemma}

\begin{proof}
For $l = 2$, if $b \equiv 1 \pmod 8$, choose $z_1 \equiv z_2 \equiv z_3 \equiv 1 \pmod 8$. If $b \equiv -1 \pmod 8$, choose $z_1 \equiv 0, z_2 \equiv z_3 \equiv 1 \pmod 8$. Hensel's lemma lifts.

For $l = p$ and $l = q$, the same arguments as in Lemma \ref{lem:exist_odd_b} apply.
\end{proof}

\begin{lemma}\label{lem:exist_even_12}
For $m$ even, if every prime factor of $n$ is $\equiv \pm 1 \pmod 8$ (i.e., $n_2 = 1$), then $(1, 2) \in \operatorname{Sel}_2(E_{m,n}/\mathbb{Q})$.
\end{lemma}

\begin{proof}
For $l \notin \{2, p, q\}$, the Hasse-Weil bound gives a solution modulo $l$ that lifts by Hensel's lemma.

For $l = 2$, write $m = 2r$. Since all prime factors of $n$ are $\equiv \pm 1 \pmod 8$, we have $q \equiv 1 \pmod 8$ (as $q \equiv 5 \pmod 8$ is impossible by quadratic reciprocity). Thus $\sqrt{q} \in \mathbb{Z}_2$.

The system reduces to finding $u, v, w \in \mathbb{Z}_2$ satisfying
\[ u^2 - 2v^2 = n^2, \quad w^2 - v^2 = q, \quad u^2 - 2w^2 = -1. \]
Since $n^2 \equiv q \equiv 1 \pmod 8$, a solution modulo $2^{r+2}$ is given by $u \equiv 1, v \equiv 0, w \equiv \sqrt{q} \pmod{2^{r+2}}$. This lifts by Hensel's lemma to $\mathbb{Z}_2$.

For $l = p$ and $l = q$, the same arguments as in Lemma \ref{lem:exist_odd_b} apply.
\end{proof}

\section{Computational Verification}
The following tables provide explicit examples computed using SageMath \cite{sage} and MAGMA \cite{magma}. The 2-Selmer ranks were computed using 2-descent routines; the generators listed are those predicted by Theorems \ref{thm:odd_m} and \ref{thm:even_m}.

Each pair $(b_1, b_2)$ represents a class in $\mathbb{Q}(S, 2) \times \mathbb{Q}(S, 2)$ modulo squares. These representatives are not rational points on the elliptic curve, but correspond to homogeneous spaces locally solvable everywhere. The 2-Selmer rank is computed as $\dim_{\mathbb{F}{2}} \operatorname{Sel}{2}(E_{m,n}/\mathbb{Q}) - 2$.

\begin{table}[htbp]
\centering
\scriptsize
\caption{Odd $m$ (representative: $m=1$): 2-Selmer ranks}
\label{tab:odd_m_ranks}
\begin{tabular}{lrr}
\toprule
$n$ & $q \pmod 8$ & $s_{2}(E_{m,n})$ \\
\midrule
3 & 5 & 0 \\
5 & 5 & 0 \\
$3 \cdot 5$ & 1 & 1 \\
$5 \cdot 7$ & 5 & 1 \\
$3 \cdot 7 \cdot 11$ & 1 & 2 \\
$5 \cdot 7 \cdot 17$ & 5 & 2 \\
$3 \cdot 7 \cdot 17 \cdot 43$ & 1 & 3 \\
$7 \cdot 17 \cdot 23 \cdot 53$ & 5 & 3 \\
$3 \cdot 7 \cdot 11 \cdot 17 \cdot 23$ & 1 & 4 \\
$7 \cdot 17 \cdot 23 \cdot 31 \cdot 53$ & 5 & 4 \\
$5 \cdot 7 \cdot 11 \cdot 17 \cdot 23 \cdot 31$ & 1 & 5 \\
$5 \cdot 7 \cdot 17 \cdot 23 \cdot 31 \cdot 41$ & 5 & 5 \\
$5 \cdot 7 \cdot 11 \cdot 17 \cdot 23 \cdot 31 \cdot 41$ & 1 & 6 \\
$7 \cdot 13 \cdot 17 \cdot 23 \cdot 31 \cdot 41 \cdot 71$ & 5 & 6 \\
$3 \cdot 5 \cdot 7 \cdot 17 \cdot 23 \cdot 31 \cdot 41 \cdot 47$ & 1 & 7 \\
$3 \cdot 7 \cdot 17 \cdot 23 \cdot 31 \cdot 41 \cdot 47 \cdot 73$ & 5 & 7 \\
$3 \cdot 7 \cdot 13 \cdot 17 \cdot 23 \cdot 31 \cdot 41 \cdot 47 \cdot 71$ & 1 & 8 \\
$7 \cdot 17 \cdot 23 \cdot 31 \cdot 41 \cdot 47 \cdot 53 \cdot 71 \cdot 73$ & 5 & 8 \\
$5 \cdot 7 \cdot 11 \cdot 17 \cdot 23 \cdot 31 \cdot 41 \cdot 47 \cdot 71 \cdot 73$ & 1 & 9 \\
$7 \cdot 17 \cdot 19 \cdot 23 \cdot 31 \cdot 41 \cdot 47 \cdot 71 \cdot 79 \cdot 89$ & 5 & 9 \\
$3 \cdot 7 \cdot 11 \cdot 17 \cdot 23 \cdot 31 \cdot 41 \cdot 47 \cdot 71 \cdot 73 \cdot 79$ & 1 & 10 \\
$7 \cdot 17 \cdot 23 \cdot 31 \cdot 37 \cdot 41 \cdot 47 \cdot 71 \cdot 73 \cdot 79 \cdot 89$ & 5 & 10 \\
$3 \cdot 5 \cdot 7 \cdot 17 \cdot 23 \cdot 31 \cdot 41 \cdot 47 \cdot 71 \cdot 73 \cdot 79 \cdot 89$ & 1 & 11 \\
$7 \cdot 17 \cdot 19 \cdot 23 \cdot 31 \cdot 41 \cdot 47 \cdot 71 \cdot 73 \cdot 79 \cdot 89 \cdot 97$ & 5 & 11 \\
\bottomrule
\end{tabular}
\end{table}

\begin{remark}
For even $m$, the 2-Selmer rank is independent of $q \pmod 8$, depending solely on $r, s$, and $t$ (the prime factor counts of $n_1$ and $n_2$). The column $q \pmod 8$ is included in Table \ref{tab:even_m_ranks} only for completeness and pattern illustration.
\end{remark}

\begin{remark}
When $n$ has $r$ distinct prime factors all congruent to $1 \pmod 8$ (so $n_2 = 1, s = 0$), the 2-Selmer rank grows linearly with $r$. For $q \equiv 1 \pmod 8$, both odd and even $m$ yield $s_2(E_{m,n}) = r + 1$, confirming that the 2-Selmer rank can be made arbitrarily large by choosing sufficiently many prime factors congruent to $1 \pmod 8$.
\end{remark}

\begin{table}[htbp]
\centering
\scriptsize
\caption{Even $m$ (representative: $m=0$): 2-Selmer ranks}
\label{tab:even_m_ranks}
\begin{tabular}{lrr}
\toprule
$n$ & $q \pmod 8$ & $s_{2}(E_{m,n})$ \\
\midrule
3 & 5 & 0 \\
5 & 5 & 0 \\
$3 \cdot 5$ & 1 & 1 \\
$3 \cdot 13$ & 1 & 1 \\
$5 \cdot 7$ & 5 & 1 \\
$79$ & 1 & 2 \\
$3 \cdot 7 \cdot 19$ & 1 & 2 \\
$3 \cdot 5 \cdot 11$ & 5 & 2 \\
$17 \cdot 23$ & 1 & 3 \\
$3 \cdot 7 \cdot 17 \cdot 43$ & 1 & 3 \\
$7 \cdot 17 \cdot 23 \cdot 53$ & 5 & 3 \\
$7 \cdot 17 \cdot 31$ & 1 & 4 \\
$3 \cdot 7 \cdot 11 \cdot 17 \cdot 23$ & 1 & 4 \\
$7 \cdot 17 \cdot 23 \cdot 31 \cdot 53$ & 5 & 4 \\
$5 \cdot 7 \cdot 11 \cdot 17 \cdot 23 \cdot 31$ & 1 & 5 \\
$5 \cdot 7 \cdot 17 \cdot 23 \cdot 31 \cdot 41$ & 5 & 5 \\
$5 \cdot 7 \cdot 11 \cdot 17 \cdot 23 \cdot 31 \cdot 41$ & 1 & 6 \\
$7 \cdot 13 \cdot 17 \cdot 23 \cdot 31 \cdot 41 \cdot 71$ & 5 & 6 \\
$3 \cdot 5 \cdot 7 \cdot 17 \cdot 23 \cdot 31 \cdot 41 \cdot 47$ & 1 & 7 \\
$3 \cdot 7 \cdot 17 \cdot 23 \cdot 31 \cdot 41 \cdot 47 \cdot 73$ & 5 & 7 \\
$3 \cdot 7 \cdot 13 \cdot 17 \cdot 23 \cdot 31 \cdot 41 \cdot 47 \cdot 71$ & 1 & 8 \\
$7 \cdot 17 \cdot 23 \cdot 31 \cdot 41 \cdot 47 \cdot 53 \cdot 71 \cdot 73$ & 5 & 8 \\
$5 \cdot 7 \cdot 11 \cdot 17 \cdot 23 \cdot 31 \cdot 41 \cdot 47 \cdot 71 \cdot 73$ & 1 & 9 \\
$7 \cdot 17 \cdot 19 \cdot 23 \cdot 31 \cdot 41 \cdot 47 \cdot 71 \cdot 79 \cdot 89$ & 5 & 9 \\
$3 \cdot 7 \cdot 11 \cdot 17 \cdot 23 \cdot 31 \cdot 41 \cdot 47 \cdot 71 \cdot 73 \cdot 79$ & 1 & 10 \\
$7 \cdot 17 \cdot 23 \cdot 31 \cdot 37 \cdot 41 \cdot 47 \cdot 71 \cdot 73 \cdot 79 \cdot 89$ & 5 & 10 \\
$3 \cdot 5 \cdot 7 \cdot 17 \cdot 23 \cdot 31 \cdot 41 \cdot 47 \cdot 71 \cdot 73 \cdot 79 \cdot 89$ & 1 & 11 \\
$7 \cdot 17 \cdot 19 \cdot 23 \cdot 31 \cdot 41 \cdot 47 \cdot 71 \cdot 73 \cdot 79 \cdot 89 \cdot 97$ & 5 & 11 \\
\bottomrule
\end{tabular}
\end{table}

\begin{table}[htbp]
\centering
\scriptsize
\caption{Odd $m$ (representative: $m=1$): Generators in $\operatorname{Sel}2(E{m,n}/\mathbb{Q})$}
\label{tab:odd_m_gens}
\begin{tabularx}{\textwidth}{lX}
\toprule
$n$ & Generators \\
\midrule
3 & $E[2]^*$ \\
5 & $E[2]^*$ \\
$3 \cdot 5$ & $(2,2)$ \\
$5 \cdot 7$ & $(14,7)$ \\
$3 \cdot 7 \cdot 11$ & $(2,2), (14,7)$ \\
$5 \cdot 7 \cdot 17$ & $(14,7), (17,17)$ \\
$3 \cdot 7 \cdot 17 \cdot 43$ & $(2,2), (14,7), (17,17)$ \\
$7 \cdot 17 \cdot 23 \cdot 53$ & $(14,7), (17,17), (46,23)$ \\
$3 \cdot 7 \cdot 11 \cdot 17 \cdot 23$ & $(2,2), (14,7), (17,17), (46,23)$ \\
$7 \cdot 17 \cdot 23 \cdot 31 \cdot 53$ & $(14,7), (17,17), (46,23), (62,31)$ \\
$5 \cdot 7 \cdot 11 \cdot 17 \cdot 23 \cdot 31$ & $(2,2), (14,7), (17,17), (46,23), (62,31)$ \\
$5 \cdot 7 \cdot 17 \cdot 23 \cdot 31 \cdot 41$ & $(14,7), (17,17), (46,23), (62,31), (41,41)$ \\
$5 \cdot 7 \cdot 11 \cdot 17 \cdot 23 \cdot 31 \cdot 41$ & $(2,2), (14,7), (17,17), (46,23), (62,31), (41,41)$ \\
$7 \cdot 13 \cdot 17 \cdot 23 \cdot 31 \cdot 41 \cdot 71$ & $(14,7), (17,17), (46,23), (62,31), (41,41), (142,71)$ \\
$3 \cdot 5 \cdot 7 \cdot 17 \cdot 23 \cdot 31 \cdot 41 \cdot 47$ & $(2,2), (14,7), (17,17), (46,23), (62,31), (41,41), (94,47)$ \\
$3 \cdot 7 \cdot 17 \cdot 23 \cdot 31 \cdot 41 \cdot 47 \cdot 73$ & $(14,7), (17,17), (46,23), (62,31), (41,41), (94,47), (73,73)$ \\
$3 \cdot 7 \cdot 13 \cdot 17 \cdot 23 \cdot 31 \cdot 41 \cdot 47 \cdot 71$ & $(2,2), (14,7), (17,17), (46,23), (62,31), (41,41), (94,47), (142,71)$ \\
$7 \cdot 17 \cdot 23 \cdot 31 \cdot 41 \cdot 47 \cdot 53 \cdot 71 \cdot 73$ & $(14,7), (17,17), (46,23), (62,31), (41,41), (94,47), (142,71), (73,73)$ \\
$5 \cdot 7 \cdot 11 \cdot 17 \cdot 23 \cdot 31 \cdot 41 \cdot 47 \cdot 71 \cdot 73$ & $(2,2), (14,7), (17,17), (46,23), (62,31), (41,41), (94,47), (142,71), (73,73)$ \\
$7 \cdot 17 \cdot 19 \cdot 23 \cdot 31 \cdot 41 \cdot 47 \cdot 71 \cdot 79 \cdot 89$ & $(14,7), (17,17), (46,23), (62,31), (41,41), (94,47), (142,71), (158,79), (89,89)$ \\
$3 \cdot 7 \cdot 11 \cdot 17 \cdot 23 \cdot 31 \cdot 41 \cdot 47 \cdot 71 \cdot 73 \cdot 79$ & $(2,2), (14,7), (17,17), (46,23), (62,31), (41,41), (94,47), (142,71), (73,73), (158,79)$ \\
$7 \cdot 17 \cdot 23 \cdot 31 \cdot 37 \cdot 41 \cdot 47 \cdot 71 \cdot 73 \cdot 79 \cdot 89$ & $(14,7), (17,17), (46,23), (62,31), (41,41), (94,47), (142,71), (73,73), (158,79), (89,89)$ \\
$3 \cdot 5 \cdot 7 \cdot 17 \cdot 23 \cdot 31 \cdot 41 \cdot 47 \cdot 71 \cdot 73 \cdot 79 \cdot 89$ & $(2,2), (14,7), (17,17), (46,23), (62,31), (41,41), (94,47), (142,71), (73,73), (158,79), (89,89)$ \\
$7 \cdot 17 \cdot 19 \cdot 23 \cdot 31 \cdot 41 \cdot 47 \cdot 71 \cdot 73 \cdot 79 \cdot 89 \cdot 97$ & $(14,7), (17,17), (46,23), (62,31), (41,41), (94,47), (142,71), (73,73), (158,79), (89,89), (97,97)$ \\
\bottomrule
\end{tabularx}
\end{table}

\begin{table}[htbp]
\centering
\scriptsize
\caption{Even $m$ (representative: $m=0$): Generators in $\operatorname{Sel}2(E{m,n}/\mathbb{Q})$}
\label{tab:even_m_gens}
\begin{tabularx}{\textwidth}{lX}
\toprule
$n$ & Generators \\
\midrule
3 & $E[2]^*$ \\
5 & $E[2]^*$ \\
$3 \cdot 5$ & $(15,15)$ \\
$3 \cdot 13$ & $(39,39)$ \\
$5 \cdot 7$ & $(7,7)$ \\
$79$ & $(79,79), (1,2)$ \\
$3 \cdot 7 \cdot 19$ & $(7,7), (57,57)$ \\
$3 \cdot 5 \cdot 11$ & $(15,15), (33,33)$ \\
$17 \cdot 23$ & $(17,17), (23,23), (1,2)$ \\
$3 \cdot 7 \cdot 17 \cdot 43$ & $(7,7), (17,17), (129,129)$ \\
$7 \cdot 17 \cdot 23 \cdot 53$ & $(7,7), (17,17), (23,23)$ \\
$7 \cdot 17 \cdot 31$ & $(7,7), (17,17), (31,31), (1,2)$ \\
$3 \cdot 7 \cdot 11 \cdot 17 \cdot 23$ & $(7,7), (17,17), (23,23), (33,33)$ \\
$7 \cdot 17 \cdot 23 \cdot 31 \cdot 53$ & $(7,7), (17,17), (23,23), (31,31)$ \\
$5 \cdot 7 \cdot 11 \cdot 17 \cdot 23 \cdot 31$ & $(7,7), (17,17), (23,23), (31,31), (55,55)$ \\
$5 \cdot 7 \cdot 17 \cdot 23 \cdot 31 \cdot 41$ & $(7,7), (17,17), (23,23), (31,31), (41,41)$ \\
$5 \cdot 7 \cdot 11 \cdot 17 \cdot 23 \cdot 31 \cdot 41$ & $(7,7), (17,17), (23,23), (31,31), (41,41), (55,55)$ \\
$7 \cdot 13 \cdot 17 \cdot 23 \cdot 31 \cdot 41 \cdot 71$ & $(7,7), (17,17), (23,23), (31,31), (41,41), (71,71)$ \\
$3 \cdot 5 \cdot 7 \cdot 17 \cdot 23 \cdot 31 \cdot 41 \cdot 47$ & $(7,7), (15,15), (17,17), (23,23), (31,31), (41,41), (47,47)$ \\
$3 \cdot 7 \cdot 17 \cdot 23 \cdot 31 \cdot 41 \cdot 47 \cdot 73$ & $(7,7), (17,17), (23,23), (31,31), (41,41), (47,47), (73,73)$ \\
$3 \cdot 7 \cdot 13 \cdot 17 \cdot 23 \cdot 31 \cdot 41 \cdot 47 \cdot 71$ & $(7,7), (17,17), (23,23), (31,31),(39,39), (41,41), (47,47), (71,71)$ \\
$7 \cdot 17 \cdot 23 \cdot 31 \cdot 41 \cdot 47 \cdot 53 \cdot 71 \cdot 73$ & $(7,7), (17,17), (23,23), (31,31), (41,41), (47,47), (71,71), (73,73)$ \\
$5 \cdot 7 \cdot 11 \cdot 17 \cdot 23 \cdot 31 \cdot 41 \cdot 47 \cdot 71 \cdot 73$ & $(7,7), (17,17), (23,23), (31,31), (41,41), (47,47), (55,55), (71,71), (73,73)$ \\
$7 \cdot 17 \cdot 19 \cdot 23 \cdot 31 \cdot 41 \cdot 47 \cdot 71 \cdot 79 \cdot 89$ & $(7,7), (17,17), (23,23), (31,31), (41,41), (47,47), (71,71), (79,79), (89,89)$ \\
$3 \cdot 7 \cdot 11 \cdot 17 \cdot 23 \cdot 31 \cdot 41 \cdot 47 \cdot 71 \cdot 73 \cdot 79$ & $(7,7), (17,17), (23,23), (31,31), (41,41), (47,47), (33,33), (71,71), (73,73), (79,79)$ \\
$7 \cdot 17 \cdot 23 \cdot 31 \cdot 37 \cdot 41 \cdot 47 \cdot 71 \cdot 73 \cdot 79 \cdot 89$ & $(7,7), (17,17), (23,23), (31,31), (41,41), (47,47), (71,71), (73,73), (79,79), (89,89)$ \\
$3 \cdot 5 \cdot 7 \cdot 17 \cdot 23 \cdot 31 \cdot 41 \cdot 47 \cdot 71 \cdot 73 \cdot 79 \cdot 89$ & $(7,7), (15,15), (17,17), (23,23), (31,31), (41,41), (47,47), (71,71), (73,73), (79,79), (89,89)$ \\
$7 \cdot 17 \cdot 19 \cdot 23 \cdot 31 \cdot 41 \cdot 47 \cdot 71 \cdot 73 \cdot 79 \cdot 89 \cdot 97$ & $(7,7), (17,17), (23,23), (31,31), (41,41), (47,47), (71,71), (73,73), (79,79), (89,89), (97,97)$ \\
\bottomrule
\end{tabularx}
\end{table}

\section*{Acknowledgments}
The first author acknowledges the support from DY Patil International University, Akurdi, Pune, India (CISR/2025SEPT/SM/007). The second author acknowledges support from SRM University, AP. The authors thank Prof. Debopam Chakraborty for his valuable suggestions during the early stages of this work.

\noindent\textit{Centre for Interdisciplinary Studies and Research, DY Patil International University, Akurdi, Pune, Maharashtra 411044, India} \\
Email: \texttt{vinodkumar.ghale@dypiu.ac.in} \\[1em]
\noindent\textit{Department of Mathematics, SRM University-AP, Amaravati 522240, Andhra Pradesh, India} \\
Email: \texttt{imdad696.ii@gmail.com}

\end{document}